\documentclass[12pt]{article}
\usepackage{graphicx}
\usepackage{amssymb}
\usepackage{verbatim}
\font\elevenss=cmss11

\font\eightss=cmss8

\font\sixss=cmss8 at 6pt

\newfam\ssfam
\textfont\ssfam=\elevenss \scriptfont\ssfam=\eightss
 \scriptscriptfont\ssfam=\sixss
\def\ss{\fam\ssfam \elevenss}%

\newtheorem {thm}{Theorem}
\newtheorem {lem}[thm]{Lemma}
\newtheorem {pr}[thm]{Proposition}

\def\Cox{\hfill \Box}
\def\sign{\mbox{sgn}}

\def\ee{\epsilon}

\def\E{{\bf{E}}}
\def\P{{\bf{P}}}
\def\R{{\mathbf{R}}}
\def\Z{{\mathbf{Z}}}

\def\F{{\cal{F}}}
\def\|{\, | \, }
\def\one{{\bf 1}}

\def\T{{\cal T}}
\def\Yt{\tilde{Y}}
\def\Mt{\tilde{M}}

\def\spce{{\cal S}}
\def\deg{{\rm deg \,}}
\def\P{{\bf P}}
\def\Q{{\bf Q}}
\def\deg{\beta}
\def\vari{\mbox{\ss span}}
\def\maxi{\mbox{\ss max}}

\begin{document}

\begin{titlepage}

\begin{center}
{\Large \bf A probabilistic model for the degree of the cancellation 
polynomial in Gosper's Algorithm}
\end{center}

\begin{flushright}
Robin Pemantle \footnote{Research supported in part by National
Science Foundation grant \# DMS 0103635}$^,$\footnote{Department
of Mathematics, The Ohio State University, 231 W. 18 Ave.,
Columbus, OH 43210}
\end{flushright}

\vfill

{\bf ABSTRACT:} \hfill \break
Milenkovic and Compton in 2002 gave an analysis of the run time of Gosper's 
algorithm applied to a random input.  The main part of this was an 
asymptotic analysis of the random degree of the cancellation polynomial 
$c(k)$ under various stipulated laws for the input.  Their methods use 
probabilistic transform techniques.  Here, a more general class
of input distributions is considered, and limit laws of the type
proved by Milenkovic and Compton are shown to follow from a 
general functional central limit theorem.  The methods herein are
probabilistic and elementary and may be used to compute the means 
of the limiting distributions. 

\vfill

 
\end{titlepage}

\section{Introduction}

Great strides have been made recently in automatic summation of
series, particularly hypergeometric series.  A source for this 
is~\cite{PWZ96}, which includes a historical development of the
problem as well as a fine exposition of the recent and seminal
work of the three authors.  A cornerstone of the automation of
hypergeometric summation is Gosper's algorithm.  In~\cite{MC1},
it is pointed out that ``despite the fact that Gosper's algorithm is 
one of the most important achievements in computer algebra, to date
there are no results concerning the average running time of the
algorithm.''  In that same work, Milenkovic and Compton undertake
an analysis of the run time under various stipulated probabilistic
models for the inputs.

To describe the results of~{MC1}, let $f$ and $g$ be polynomials,
let $r_k := f(k) / g(k)$, and let $t_n := \prod_{k=1}^{n-1} r_k$.
The series $\{ t_n \}$ and its partial sums are known as {\em 
hypergeometric}, and $r_k$ is called the {\em hypergeometric ratio}.
The purpose of Gosper's algorithm is to find a closed form expression
for the partial sum $S_N := \sum_{n=1}^N t_n$.  Its input is often 
specified as the rational function $r_k$ in factored form.  Roots
of $f$ and $g$ differing by integers play a crucial role in the 
algorithm.  Milenkovic and Compton observe that not much generality 
is lost in assuming the roots have been classified according to 
their remainders modulo~1, and that the problem has been restricted 
to one of these moduli classes.  In other words, they assume that 
$f$ and $g$ have integer roots.  They go on to stipulate joint probability
distributions for $f$ and $g$, which have as parameters an {\em a priori}
bound on the location of the roots.  Specifically, they assume that
\begin{eqnarray}
f(k) & = & \prod_{j=1}^m (k-j)^{A_j} \, \nonumber \\[1ex]
g(k) & = & \prod_{j=1}^m (k-j)^{B_j} \, \label{eq:fg}
\end{eqnarray}
so all roots of $f$ and $g$ lie in $[m] := \{ 1 , \ldots , m \}$. 

The ``uniform {\bf R} model'' considered by Milenkovic and Compton
may be described as follows.  For each of $f$ and $g$, a sequence of
$n$ IID uniform picks is made from $[m]$.  This gives the multisets of
roots for $f$ and $g$.  In other words, for $j \in [m]$, the random
variables $A_j$ and $B_j$ count how many of $n$ IID picks from $[m]$
are equal to $j$.  Milenkovic and Compton point out that this is an 
urn model of Maxwell-Boltzman type.  A key to their analysis is the
representation of the variables $\{ A_j \}$ as a set of $m$ IID picks
from a Poisson distribution (of any mean), conditioned to sum to $n$.  
A second distribution of roots they analyze is the ``multi-set {\bf R} 
model'', in which the multisets of roots (equivalently the sequences 
$(A_j)_{j=1}^m$ and $(B_j)_{j=1}^m$) are chosen uniformly from all 
multisets (equivalently all sequences of $m$ nonnegative integers 
summing to $n$).  This is a Bose-Einstein urn model, 
and is equivalent to conditioning two IID sequences of geometric 
random variables (with any mean) both to sum to $n$.  They also
discuss two models, the ``Uniform {\bf T} model'' and the ``Multiset
{\bf T} model'', in which the roots of the numerator and denominator
of the partial product terms $t_k$ are directly modeled by the
two respective urn models; these models are not addressed in this
paper.  

Milenkovic and Compton give a partial average case analysis, meaning 
that they focus on a few quantities which are highly determinative
of the run time and give average case analyses of these.  The most
important such quantity is the degree of the {\em cancellation polynomial},
$c(k)$.  This is defined as the minimal polynomial $c$ for which we 
may write
\begin{equation} \label{eq:cancellation}
\frac{f(k)}{g(k)} = \frac{a(k)}{b(k)} \frac{c(k+1)}{c(k)}
\end{equation}
and also satisfying
\begin{equation} \label{eq:no integer}
{\rm GCD} (a(x) , b(x-h)) = 1 \mbox{ for all nonnegative integers} h \, .
\end{equation}
The determination of this 
polynomial is Step~2 in the version of Gosper's algorithm described
in~\cite{Wis03}, which is distilled from~\cite{PWZ96}.  
Milenkovic and Compton obtain the following results (they use $N$ in 
place of the $m$ in this paper).  The draft of their manuscript 
cited here is a very preliminary version which the authors have kindly 
provided.  Consequently, only results independently proved in the present
paper are quoted here, though in fact the manuscript~\cite{MC1}
obtains explicit constants for the asymptotic expectations.
\begin{thm}[Milenkovic and Compton (2002) Theorems~20 and~21] \label{th:MC1}
In the \hfill \\ uniform {\bf R} model, if $n/m \to \lambda$, 
the expected degree of $c(k)$ is asymptotic to
$$C_1 (\lambda) m^{3/2} \, .$$
In the multiset {\bf R} model, they find that when $\lambda$ is 
sufficiently large, the expected degree of $c$ is asymptotically
$$ C_2 (\lambda) m^{3/2} \, .$$
\end{thm}

The method of~\cite{MC1} is to compute transforms (generating
functions) for the unconditioned distributions, in which the variables
$A_j$ and $B_j$ are independent Poissons or geometrics, and then
de-Poissonize, according to machinery they developed in~\cite{MC2}.

Analytic de-poissonization may be technically somewhat involved; see
for example~\cite{JS99}.  The present paper also relies on the 
representation of the stipulated distributions as IID conditioned 
on a fixed sum.  After that, however, the method herein is purely 
probabilistic, relying on limit theory for the random walk whose
increments are $A_j - B_j$.  Theorem~\ref{th:soft} below, whose proof
is a straightforward application of random walk limit theory,
shows that as $n,m \to \infty$ with $n/m = \lambda + o(\lambda^{-1/2})$, 
the expected degree of $c$ is 
\begin{equation} \label{eq:form}
(c + o(1)) m^{3/2}
\end{equation}
where $c$ is a certain expectation taken with respect to the 
Brownian bridge.  For the uniform and multiset {\bf R} models, 
the constant $c$ is calculated respectively as
\begin{eqnarray}
c_{\rm unif} & = & \frac{\pi \sqrt{2 \lambda}}{16} \label{eq:uniform R}  
   \, ; \\[1ex]
c_{\rm multi} & = & \frac{\pi \sqrt{2 \lambda (\lambda+1)}}{16} 
   \label{eq:multiset R} \, .
\end{eqnarray}
The authors of~\cite{MC1} are aware of the random walk representation, 
but it appears that they use this only via analytic transforms,
and not via any scaling limits of the random walk paths.

\section{Definitions and results}

Let $F$ be a distribution on the nonnegative integers with mean
$\lambda$ and variance $V < \infty$.  We assume throughout that
the GCD of the support of $F$ is~1, as is true for Poisson 
distributions and geometric distributions of any mean.  
Let $\P_{F,m}$ be the probability measure on $\spce := \Z^{2m}$
making the coordinates IID with common distribution $F$.  Denote
the first $m$ coordinates by $A_1 , \ldots , A_m$ and the
last $m$ coordinates by $B_1 , \ldots , B_m$.  Let $\Q_{F,m;n}$
be the result of conditioning $\P_{F,m}$ so that
$$\sum_{j=1}^m A_j = \sum_{j=1}^m B_j = n \, .$$
Associated to each $\omega \in \spce$ are the polynomials
$f$ and $g$ defined by~(\ref{eq:fg}), and the associated
hypergeometric series with ratio $r_k = f(k) / g(k)$.   
Define a function $\deg : \spce \to \Z^+$ by letting 
$\deg (\omega)$ be the minimal degree of a polynomial $c$
satisfying~(\ref{eq:cancellation}) for the polynomials $f$ and $g$.
The next result, proved at the end of the section, provides 
an alternative expression for $\deg$, which is essentially the 
random walk representation in~\cite[Theorem~1]{MC1}.

Define $X_j := A_j - B_j$ and $S_k := \sum_{j=1}^k X_j$, with the
convention that $S_0 := 0$.  Define $M_k := \min \{ S_j : 0 \leq j
\leq k\}$ and define $\Mt_k := \min \{ S_j : k \leq j \leq m \}$.
Let $\tau := \min \{ k : S_j \geq S_k \; \forall j > k \}$
denote the time of the first minimum of the process $\{ S_j : 0 \leq j
\leq m \}$.  For each $j$, define $Y_j := S_j - M_j$ and
$\Yt_j := S_j - \Mt_j$.  Define $I_j := Y_j \one_{j \leq \tau}
+ \Yt_j \one_{j > \tau}$.  At the end of this section we will prove:

\begin{lem} \label{lem:D}
$$\deg = \sum_{j=1}^{m-1} I_j = \sum_{j=1}^\tau Y_j + \sum_{j=\tau+1}^{m-1}
   \Yt_j \, .$$
\end{lem}

Under the measures $\P_{F,m}$ and $\Q_{F,m;n}$, the quantities 
$A_j, B_j, X_j, Y_j, \deg$ and so forth become random variables.
In order to state the main results of this paper, some definitions
are required that mirror the definitions of these quantities
but on the space of continuous limits.  

Let $\Omega$ be the space of CADLAG paths on $[0,1]$, with
filtration $\{ \F_t \}$.  For $\omega \in \Omega$, define
the minimum process $M^*$ by $M^* (\omega) (t) = \inf \{ \omega (s) :
s \leq t \}$.  Define $\tilde{M^*}$ to be the right to left minimum
process $\tilde{M^*} (\omega) (t) = \inf \{ \omega (s) :
s \geq t \}$.  Define the process $Y$ to be $\omega - M^*$ and
$\tilde{Y}$ to be $\omega - \tilde{M^*}$.  Define $\tau = \inf \{ t :
\omega(s) \geq \omega (t) \, \forall s > t \}$ to be the time
of the first minimum and let $I$ be the process $Y \one_{t \leq \tau}
+ \tilde{Y} \one_{t > \tau}$.  Finally, we define $\deg^* (\omega)$
to be the quantity $\int_0^1 I(\omega) (t) \, dt$.

Let $\P_{2V}$ be the law of a centered continuous Gaussian process 
with covariances 
$$\E \omega (s) \omega (t) = 2 V (s \wedge t) \, .$$
In other words, $\P_{2V}$ is a Brownian motion with amplitude $\sqrt{2V}$.
Let $\Q_{2V}$ be the law of a centered continuous Gaussian process
with covariances
\begin{equation} \label{eq:bridge covariances}
\E \omega (s) \omega (t) = 2V ((s \wedge t) - st) \, .
\end{equation}
In other words, $\Q_{2V}$ is a Brownian bridge of amplitude $\sqrt{2V}$. 
Recall that the mean and variance of $F$ are denoted by $\lambda$ and $V$.
The first main result of the paper, proved in the next section, is:
\begin{thm} \label{th:soft}
Let $m \to \infty$ with $n = \lambda m + o(m^{1/2})$.  Then
\begin{equation} 
\int \deg \, d\P_{F,m} = (K_1 + o(1)) m^{3/2} 
   \label{eq:soft1} \\[1ex]
\end{equation}
and
\begin{equation}
\int \deg \, d\Q_{F,m;n} = (K_2 + o(1)) 
   m^{3/2} \label{eq:soft2}
\end{equation}
with
\begin{eqnarray}
K_1 := \int \deg^* \, d \P_{2V} \, ; \label{eq:K_1} \\[1ex]
K_2 := \int \deg^* \, d \Q_{2V} \, . \label{eq:K_2} 
\end{eqnarray}
\end{thm}

The integrals~(\ref{eq:K_1}) and~(\ref{eq:K_2}) may be evaluated,
leading to quantitative versions:

\begin{thm} \label{th:IID}
Consider the measure $\P_{F,m}$, under which
$\{ A_j , B_j : 1 \leq j \leq m \}$ are IID and have common variance
$V$.  Then as $m \to \infty$,
$$\E \deg = (1 + o(1)) \frac{2}{3} \sqrt{\frac{V}{\pi}} m^{3/2} \, .$$
\end{thm}

\begin{thm} \label{th:cond}
Suppose that $N(m)$ is an integer satisfying $n(m) / m = \lambda + o(m^{-1/2})$ 
as $m \to \infty$.  Consider the measure $\Q_{F,m;n}$, under which
$\{ A_j , B_j : 1 \leq j \leq m \}$ have the distribution of IID picks 
from $F$ conditioned on $\sum_{j=1}^m A_j = \sum_{j=1}^m B_j = N$.  
Then as $m \to \infty$,
$$\E \deg = (1 + o(1)) \frac{\pi \sqrt{2}}{16} \sqrt{V} m^{3/2} \, .$$
\end{thm}

\noindent{\em Remarks:} 
\begin{quote}
(1) Theorem~\ref{th:soft} emerges without much difficulty from the 
convergence of the random walks to Brownian paths.  Thus not only 
does $m^{-3/2}$ times the $\P_{F,m}$ expectation of $\deg$ converge 
to $\E_{2V} \deg^*$, but the $\P_{F,m}$ distribution of $m^{-3/2} \deg$ 
converges to the $\P_{2V}$ distribution of $\deg^*$. \hfill \\

(2) Although it is relatively easy and is superseded by the two 
quantitative results, Theorem~\ref{th:soft} is worth stating separately
for the following reason:
the computations in Theorems~\ref{th:IID} and~\ref{th:cond} are 
a little tricky, and it is instructive to see that the form of the 
result does not depend on calculations which are not transparent.  \hfill \\

(3) The asymptotics obtained by~\cite{MC1} for both the uniform and
multiset models pertain to Theorem~\ref{th:cond}.  Theorem~\ref{th:IID}
corresponds to two models discussed in~\cite{MC1} but not quoted here,
where the roots of $f$ and $g$ are assumed to be unconditioned picks
from the Poisson (respectively geometric) distributions, hence not 
necessarily equinumerous.
\end{quote}

The rest of the organization of this paper is as follows.
The functional central limit arguments are spelled out in 
Section~\ref{ss:soft}.  The values of the constants $K_1$ and $K_2$ 
are then computed in Section~\ref{ss:constants}.
In the remainder of this section we prove Lemma~\ref{lem:D}.  We
begin with an intermediate representation.

For each $j$ between 1 and $m$, place $A_j$ red balls and $B_j$ blue
balls in an urn marked $j$; this will be called ``position $j$''
or ``time $j$''.  An {\em admissible} matching of the balls is a set
of pairs of balls such that
\begin{enumerate}
\item each pair contains one red ball and one blue ball; say that these
two balls are {\em matched}, and any ball not in the union of pairs is
called {\em unmatched};
\item the pairs are pairwise disjoint;
\item if a red ball in position $i$ is matched with a blue ball in
position $j$, then $i \leq j$;
\item if there is an unmatched red ball in position $i$ and an unmatched
blue ball in position $j$ then $i > j$ (that is, among unmatched balls,
all red balls sit to the right of all blue balls).
\end{enumerate}
For each admissible matching $\xi$, define the weight $w(\xi)$ to be
the sum over all pairs in the matching of $j-i$, where $j$ is the
position of the blue ball and $i$ is the position of the red ball.

\begin{pr} \label{pr:deg}
$$\deg = \min \{ w(\xi) : \xi \mbox{ is an admissible matching} \}.$$
\end{pr}

\noindent{\sc Proof:} For every matched pair in positions $i$ and $j$,
we have
$$\frac{x-i}{x-j} = \prod_{s=i}^{j-1} \frac{x-s}{x-(s+1)} = 
   \frac{e(k+1)}{e(k)}$$
where $e(x) = \prod_{s=i}^{j-1} (x-s)$.  Thus any admissible matching 
$\xi$ of weight $w$ yields a solution $H(\xi)$ to~(\ref{eq:cancellation}) 
where $c$ has degree $w$.  Conversely, let $c$ solve~(\ref{eq:cancellation}) 
and have degree $w$.  If $w=0$ then the empty matching is admissible.  
Assume now, for induction, that $w > 0$ and that for any $f',g'$ represented
by red and blue balls in urns, and any solution $c'$ 
to~(\ref{eq:cancellation}) of degree less than $w$, there
is an admissible $\xi$ for that urn problem with $c' = H(\xi)$.  
For some $j$ which is a root of $c$, let $c'$ denote $c$ with 
the linear factor $(x-j)$ removed.  By induction, there is a 
matching $\xi'$ of weight $w-1$, with $H(\xi') = c'$, admissible for
the urn problem gotten from the original one by adding a blue ball 
at position $j$ and a red ball at $j+1$.  Both of these new balls
must be matched in $\xi'$, say in pairs of positions $(s,j)$ and $(j+1,t)$,
whence replacing these with the pair $(s,t)$ produces a $\xi$ admissible
for the original problem with $H(\xi) = c$.  $\Cox$

\noindent{\sc Proof of Lemma}~\ref{lem:D}:
For any admissible matching $\xi$, define $d_j (\xi)$
to be the number of red balls in positions $1 , \ldots , j$ matched with
blue balls in positions $j+1 , \ldots , m$.  Elementarily, summing by parts,
\begin{equation} \label{eq:w}
w(\xi) = \sum_{\parbox{0.8in}{\scriptsize pairs of \\[-1ex] positions 
   $\; (i,j)$}} \hspace{.3in} \sum_{j \leq t < i} 1 = \sum_t d_t (\xi) \, .
\end{equation}
Fix any $j \leq \tau$ and let $j_* = {\rm argmin} \{ S_j : 0 \leq i
\leq j \}$ so that $S_{j_*} = M_j$.  Then in positions
$j_* + 1 , \ldots , j$ there are a total of $Y_j$ more red balls
than blue balls.  In positions $j+1 , \ldots , \tau$, there is an excess of
$S_j - S_\tau \geq S_j - M_j = Y_j$ blue balls over red balls.
In an admissible matching, either every red ball at a position at
most $j$ is matched or every blue ball at a position at least $j+1$
is matched.  It follows that either at least $Y_j$ red balls from
positions $j_* + 1 , \ldots j$ are matched to blue balls in positions
beyond $j$, or at least $Y_j$ blue balls from positions $j+1 , \ldots , \tau$
are matched with red ball in positions up to $j$.  In either case,
$d_j (\xi) \geq Y_j$.

Similarly, fix any $j > \tau$ and let $j^* \geq j$ satisfying $S_{j^*} =
\Mt_j$.  There is an excess of $\Yt_{j-1}$ blue balls in positions
$j , \ldots , j^* - 1$.  There is an excess of red balls in positions
$\tau + 1 , \ldots , j-1$ of $S_{j-1} - S_\tau \geq \Yt_{j-1}$.
Reasoning as before, one sees that $d_{j-1} (\xi) \geq \Yt_{j-1}$.
Summing over $j$ now gives
\begin{equation} \label{eq:Y}
\sum_{j = 1}^{m-1} d_j (m) \geq \sum_{j=1}^\tau Y_j + \sum_{j = \tau}^{m-1}
   \Yt_j \, .
\end{equation}
Note that $Y_\tau = \Yt_\tau = 0$, so this is equal to $\sum_{j=1}^{m-1}
I_j$.  Minimizing over $\xi$ then gives half of the conclusion of the
theorem: $\deg \geq \sum_{j=1}^{m-1} I_j$.

To prove the other half, we produce an admissible matching $\xi$ with
$d_j (\xi) = I_j$ for all $1 \leq j \leq m-1$.  In particular, $d_\tau = 0$,
so no ball in a position at most $j$ is matched with a ball in a position
beyond $j$, and $\xi$ may be decomposed into a matching on balls in urns
$1 , \ldots , j$ and another on balls in urns $j+1 , \ldots , m$.  We
construct these separately.  An algorithm for the first is as follows.
\begin{quote}
Initialize $j := 1$.  Initialize a LIFO stack.  Pull red balls out of
urn $j$ and place them on the stack until there are no more red balls
in urn $j$.  Pull blue balls out of urn $j$: while the stack is non-empty,
match each blue ball with the top element of the stack; once the stack
is empty, label each new blue ball ``unmatched'' and discard it.  When
urn $j$ is empty, increment $j$ and execute the loop until finished with
the step $j = \tau$.
\end{quote}

It is easy to see that all red balls in positions up to $\tau$ will
be matched, since every time a red ball goes on the stack, there are
more blues and reds to follow by time $\tau$ and no blue will be discarded
until that red ball is matched.  Inductively, it is easy to check that:
\begin{itemize}
\item the stack size after step $j$ is $Y_j$;
\item the change in stack size from time $j-1$ to $j$ is
$\max \{ X_j , -Y_j \}$;
\item the total number of balls discarded through time $j$ is $-M_j$.
\end{itemize}
From this one sees that $d_j (\xi)$ is equal to the stack size after time $j$,
and is therefore equal to $Y_j$.  A stack algorithm dual to this
works in the case $j > \tau$, working backward from time $m$ to $\tau$,
stacking blue balls and matching or discarding reds.  It constructs
the other half of $\xi$ so that $d_j (\xi) \geq Y_j$ for all $j >\tau$.
This completes the proof of Lemma~\ref{lem:D}.  $\Cox$

\section{Proof of Theorem~{\protect{\ref{th:soft}}}}
\label{ss:soft}

Let $\T$ be the topology on $\Omega$ generated by the sup norm
$|\omega| := \sup_{0 \leq t \leq 1} |\omega (t)|$.  The following
lemma is necessary only because $\tau$ is not a continuous function.

\begin{lem} \label{lem:continuous}
$\deg^*$ is a continuous function on $\Omega$ with respect to $\T$.  
\end{lem}

\noindent{\sc Proof:} Suppose $\sup_s |\omega_1 (s) - \omega_2 (s))|
\leq \ee$.  Fix any $t$.  It is clear that $|M (t , \omega_1)
- M (t , \omega_2)| \leq \ee$ and likewise for $\Mt$, hence
these are continuous.  If $\tau (\omega_1) , \tau (\omega_2) \geq t$, 
it follows immediately also that $|I(t , \omega_1) - I(t , \omega_2)| 
\leq 2 \ee$.  Likewise, if $\tau (\omega_1) , \tau (\omega_2) > t$, 
it follows that $|I(t , \omega_1) - I(t , \omega_2)| \leq 2 \ee$.

Suppose now that $\tau_1 := \tau (\omega_1) < t \leq \tau_2 :=
\tau (\omega_2)$.  Then
$$\Mt (t , \omega_1) \geq \omega_1 (\tau_1) \geq \omega_2 (\tau_1) - \ee
   \geq M(t , \omega_2) - \ee$$
since $\tau_1$ is one of the times over which the inf defining $M(t ,
\omega_2)$ is taken.  Similarly,
$$\Mt (t , \omega_2) \geq \omega_2 (\tau_2) \geq \omega_1 (\tau_2) - \ee
   \geq M(t , \omega_1) - \ee \, .$$
It follows that
\begin{equation} \label{eq:M}
\left | \Mt (t , \omega_1) - M(t , \omega_2) \right | \leq \ee \, .
\end{equation}
Together with $|\omega_1 (t) - \omega_2(t)| \leq \ee$,
this shows that $|I(t , \omega_1) - I(t , \omega_2)| \leq 2 \ee$.
A similar argument shows this in the case that $\tau_2 < t \leq\tau_1$.
This establishes continuity of $I$, with continuity of $D$ following 
by integration.   $\Cox$

\noindent{\sc Proof of~(\protect{\ref{eq:soft1})}}:
Recall the definition of the partial sums $S_j$ on $\spce$
and define a map $\kappa: \spce \to \Omega$ by
$$\kappa (\omega) (t) = m^{-1/2} S_{\lfloor m t \rfloor} \, .$$
The following relations are evident:
\begin{eqnarray*}
\tau^* \circ \kappa & = & m^{-1} \tau  \, ; \\
M^* (\kappa (\omega) (t)) & = & m^{-1/2} M(\omega (\lfloor m t \rfloor))
   \, ; \\
\tilde{M^*} (\kappa (\omega) (t)) & = & m^{-1/2} \tilde{M}
   (\omega (\lfloor m t \rfloor)) \, ; \\
Y^* (\kappa (\omega) (t)) & = & m^{-1/2} Y(\omega (\lfloor m t \rfloor))
   \, ; \\
\tilde{Y^*} (\kappa (\omega) (t)) & = & m^{-1/2} \tilde{Y}
   (\omega (\lfloor m t \rfloor)) \, ; \\
I^* (\kappa (\omega) (t)) & = & m^{-1/2} I(\omega (\lfloor m t \rfloor))
   \, ; \\
\deg^* \circ \kappa & = & m^{-3/2} \deg \, .
\end{eqnarray*}
Let $\P_{(m)}$ denote the image under $\kappa$ of $\P_{F,m}$.
The functional central limit theorem says that the laws under $\P_{(m)}$ 
of $\kappa$ converge weakly as $m \to \infty$ to the measure $\P_{2V}$. 
See~\cite[Theorem~37.8]{Bil86} for a proof when $p \geq 4$
or~\cite[Theorem~6.3 of Chapter~7]{Dur96} for a general proof
using Skorohod embedding.  This and Lemma~\ref{lem:continuous}
would complete the proof of~(\ref{eq:soft1}) if $\deg^*$ were
bounded.  Since $\deg^*$ is not bounded, we may define for each $L > 0$,
$$I_L := \sign (I^*) \left ( |I^*| \wedge L \right ) \, .$$
We may then conclude that the expectation with respect to $\P_{(m)}$
of the bounded continuous function $\deg_{(L)} := \int I_{(L)}$  
converges as $m \to \infty$ to its expectation with respect to $\P_{2V}$.

\begin{lem} \label{lem:uniform bound}
Let $\maxi := \sup_{0 \leq t \leq 1} \omega (t)$.  Then 
$$\lim_{L \to \infty} \sup_m \E_{(m)} \maxi \one_{\maxi > L} = 0 \, .$$
\end{lem}

Assuming this for the moment, we observe that $\P_{(m)}$ is symmetric
so the same holds with $\inf$ in place of $\sup$, and thus
\begin{equation} \label{eq:vari}
\lim_{L \to \infty} \sup_m \E_{(m)} \vari \one_{\vari > L} = 0
\end{equation}
where $\vari = \sup_t \omega (t) - \inf_t \omega (t)$.  Since 
$$|\deg^* (\omega) - \deg_{(L)}(\omega)| \leq \vari := \sup_t \omega (t) 
   - \inf_t \omega (t) $$
and $\deg^* = \deg_{(L)}$ on the event $\{ \vari \leq L \}$,
the inequality~(\ref{eq:soft1}) follows from~(\ref{eq:vari}) and
the convergence of $\E_{(m)} \deg_{P(L)}$ to $\E_{2V} \deg_{(L)}$.

\noindent{\sc Proof of Lemma}~\ref{lem:uniform bound}:
By the $L^2$ maximum inequality (\cite[Theorem 4.4.3]{Dur96}),
\begin{eqnarray}
\E_{(m)} \maxi^2 & \leq & 4 \E \omega(1)^2 \nonumber \\
& = & 4 \E_{F,m} (m^{-1/2} S_m)^2 \nonumber \\
& = & 4 V \label{eq:bd} \, ,
\end{eqnarray}
and hence 
\begin{equation} \label{eq:L2}
\E_{(m)} \maxi \one_{\maxi > L} \leq L^{-1} \E_{(m)} \maxi^2
   \leq \frac{4V}{L}
\end{equation}
for all $m$, proving the lemma.   $\Cox$

\noindent{\sc Proof of~(\protect{\ref{eq:soft2})}}:
Recall that $\Q_{2V}$ denotes the law on $\Omega$ of a Brownian bridge of
amplitude $\sqrt{2V}$.
The proof proceeds analogously to the proof of~(\ref{eq:soft1}).
In place of the standard functional central limit theorem is a
well known result that may be found, among other places, 
in~\cite[equation~(6) of Section~0.4]{Pit02} (refer to~\cite{DK63}
for the proof).
\begin{lem}[Conditional Functional CLT] 
Let $n(m) \to \infty$ as $m \to \infty$ with $n(m)/m = \lambda 
+ o(m^{-1/2})$.  Recall that $\Q_{F,m;n}$ is the measure on $(\Z^+)^{2m}$ 
whose coordinates have the distribution of IID draws from $F$ 
conditioned on $\sum_{j=1}^m A_j = \sum_{j=1}^m B_j = n$ and recall 
the aperiodicity assumption on $F$.  Let $\Q_{(m)}$ denote the image 
under $\kappa$ of $\Q_{F,m;n}$.  Then the $\Q_{(m)}$ law of $\{ S_j : 
1 \leq j \leq m \}$ converges weakly to $\Q_{2V}$.   $\Cox$
\end{lem}

All that remains is to show the analogue of Lemma~\ref{lem:uniform bound} 
with $\Q_{(m)}$ in place of $\P_{(m)}$.  As in~(\ref{eq:L2}), 
this will follow once we have shown a uniform bound on the $\Q_{(m)}$ 
second moment of $\maxi^2$ analogous to~(\ref{eq:bd}).  This in
turn follows immediately from the inequality
\begin{equation} \label{eq:compare}
\Q_{F,m;n} (\maxi > L) \leq C \P_{F,m} (\maxi > L)
\end{equation}
where hereafter $C$ may change from equation to equation
but will depend only on $F$.  To prove~(\ref{eq:compare}),
let $G$ be the event that $\sum_{j=1}^m A_j = \sum_{j=1}^m B_j = n$.
By the local central limit theorem, $\P_{F,m} (G)$
is asymptotic to $C / n$, hence, using time-reversal symmetry 
of the path $\{ S_j \}$ under $\P_{F,m}$,
\begin{eqnarray}
\Q_{F,m;n} (\maxi > L) & = & \frac{\P_{F,m} (\maxi > L ; G)}{\P_{F,m} (G)}
   \nonumber \\[1ex]
& \leq & C \, n \P_{F,m} (\maxi > L ; G) \nonumber \\[1ex]
& = & 2 C \, n \sum_{j > L} \sum_{t \leq m/2} 
   \P_{F,m} (\maxi = j = S_t ; G) \, . \label{eq:22}
\end{eqnarray}
It then suffices to show 
\begin{equation} \label{eq:33}
\P_{F,m} (\maxi = j = S_t ; G) \leq C n^{-1} \P_{F,m} (\maxi = j = S_t)
\end{equation}
since then resumming~(\ref{eq:22}) proves~(\ref{eq:compare}).

Let $l = \lfloor 3m/4 \rfloor$, let $\F_l$ be the $\sigma$-field generated
by $\{ A_i , B_i : i \leq l \}$ and let $H$ be the event that
$S_t = \max \{ S_i : i \leq l \}$.  The local central limit 
theorem~\cite[(5.2)~in Chapter 2]{Dur96} gives
$\P_{F,m} (G \| \F_l) \leq C n^{-1}$, whence
\begin{eqnarray*}
\P_{F,m} (\maxi = j = S_t ; G) & \leq & \P_{F,m} (S_t = j ; H \cap G) \\
& \leq & \P_{F,m} (S_t = j ; H) \P_{F,m} (G \| \F_l) \\
& \leq & C n^{-1} \P_{F,m} (S_t = j ; H) 
\end{eqnarray*}
Conditioning again on $\F_l$, we see that $\P_{F,m} (S_t = j ; H) \leq
C n^{-1} \P_{F,m} (S_t = j = \maxi)$, which establishes~(\ref{eq:33}), 
hence~(\ref{eq:compare}) and the theorem.
$\Cox$

\section{Evaluation of the constants} \label{ss:constants}

Let
$$D_1 (\omega) := \int_0^\tau (\omega (t) - M(t)) \, dt \, .$$
By symmetry, $\int D_1 \, d\P_{2V} = (1/2) \int \deg^* \, d\P_{2V}$ and
similarly, $\int D_1 \, d\Q_{2V} = (1/2) \int \deg^* \, d\Q_{2V}$.  Thus
it suffices to compute expectations of $D_1$.

\noindent{\sc Proof of Theorem}~\ref{th:IID}:
For a process with law $\P_{2V}$, the process $Y(t) = \omega (t) -
M(t)$ is well known to have the same distribution as the law under
$\P_{2V}$ of the reflected Brownian motion $\{ |\omega (t)| :
0 \leq t \leq 1 \}$ (see, e.g.,~\cite[Prop.~13.13]{Kal02}).  
Clearly, the map 
$\omega \mapsto \omega - M$ has the property that $t$ is a left-to-right
minimum for $\omega$ if and only if $t$ is a zero of $\omega - M$.
The last left-to-right minimum is the global minimum, which occurs
at the last zero of $\omega - M$.  Because the process $\omega -
M$ is distributed as $|\omega|$, the location of the last
zero of $\omega - M$ is distributed as the last zero of $\omega$.
This has an arc-sine density $\pi^{-1} (x(1-x))^{-1/2}
\, dx$~\cite[Example~4.4]{Dur96}.  The following lemma writing Brownian
motion as a mixture of bridges up to the last zero follows directly
from the strong Markov property, scaling, and the fact that the
bridge is a Brownian motion conditioned to return to zero:

\begin{lem} \label{lem:RCP}
Let $L = L (\omega) = \sup \{t \leq 1 : \omega (t) = 0 \}$ be the last 
zero of Brownian motion and let $g : \Omega \to \R^+$ depend only on 
$\omega |_{[0,L(\omega)]}$.  Then
$$\int_\Omega g \, d \P_{2V} = \int_0^1 \frac{dt}{\pi \sqrt{t(1-t)}}
   \int_\Omega g \, d\Q_{2V}^{(t)} (\omega)$$
where the inner integral is on $\F_t$ and $\Q_{2V}^{(t)}$ is the 
law on $\F_t$ of a bridge of amplitude $\sqrt{2V}$ on $[0,t]$,
that is, a centered Gaussian process with covariance 
$$\E \omega (u) \omega (v) = u \wedge v - \frac{uv}{t} \, .$$
$\Cox$
\end{lem}

Using this, $c$ may be evaluated as follows.  By definition, 
by~\cite[Prop.~13.13]{Kal02}, and lastly by Lemma~\ref{lem:RCP}
applied to the integral up to $L$ of $|\omega|$, we have
\begin{eqnarray*}
\int_\Omega D_1 \, d \P_{2V} (\omega) & = & \int_\Omega \int_0^\tau
   (\omega (s) - M(\omega)(s)) \, ds \, d\P_{2V} \\[1ex]
& = & \int_\Omega \int_0^{L(\omega)} |\omega (s)| \, ds \, d\P_{2V} \\[1ex]
& = & \int_0^1 \frac{dt}{\pi \sqrt{t(1-t)}} \int_\Omega d\Q_{2V}^{(t)} 
   \int_0^t |\omega (s)| \, ds \, .
\end{eqnarray*}
 From the covariance structure of $\Q_{2V}^{(t)}$, we see this law makes 
$\omega(s)$ is a centered Gaussian with variance $s(1-s/t)$.  The 
expected absolute value of a $N(0,a)$ random variable is $\sqrt{2a/\pi}$.
Switching the order of the two inner integrals, we may then write
\begin{eqnarray*}
\int_\Omega D_1 \, d \P_{2V} (\omega) & = & 
   \int_0^1 \frac{dt}{\pi \sqrt{t(1-t)}} \int_0^t \, ds \, 
   \int_\Omega |\omega (s)| \, d\Q_{2V}^{(t)} (\omega)  \\[1ex]
& = & \int_0^1 \frac{dt}{\pi \sqrt{t(1-t)}} \int_0^t 
   \frac{ds}{\pi \sqrt{t(1-t)}} \sqrt{\frac{4V}{\pi}} \sqrt{s
   \left ( 1-\frac{s}{t} \right ) } \, .
\end{eqnarray*}
The evaluation is now straightforward integration.  Substitute $s = tu$
and $ds = t \, du$ to get
\begin{eqnarray*}
\int_\Omega D_1 \, d \P_{2V} (\omega) & = & 
   \sqrt{\frac{4V}{\pi}} \int_0^1 \frac{dt}{\pi \sqrt{1-t}} \, \int_0^1 
   \sqrt{u(1-u)} \, du \\[1ex]
& = & \sqrt{\frac{4V}{\pi}} \int_0^1 \frac{dt}{\pi \sqrt{1-t}} 
   \frac{\pi}{8} \\[1ex]
& = & \frac{1}{4} \, \sqrt{\frac{V}{\pi}} \, \int_0^1 \frac{t}{\sqrt{1-t}} 
   \, dt \\[1ex]
& = & \frac{1}{3} \sqrt{\frac{V}{\pi}} \, .
\end{eqnarray*}

Doubling yields $\int \deg^* \, d\P_{2V}$ and finishes the proof of
Theorem~\ref{th:IID}, that is, $K_1 = (2/3) \sqrt{V/\pi} \approx 0.376
\sqrt{V}$.  $\Cox$

\noindent{\sc Proof of Theorem}~\ref{th:cond}:
If the distribution of $\omega (t) - m(t)$ for a Brownian bridge were 
explicitly known in a usable form, the computation of the $\Q$-expectation 
would be analogous to the $\P$-expectation of $D$.  In the absence
of such a representation, the second computation ignores the
representation of the law of $\omega - M$ as that of $|\omega|$
and proceeds as follows.  

The counterpart to Lemma~\ref{lem:RCP} is the following joint density
for the pair $(\omega (t) , M(\omega )(t))$ under $\Q_1$ conditioned
on $\tau > t$. 
\begin{lem} \label{lem:density}
For fixed $t \in (0,1)$, define the positive function $f_H$ on the set
$R := \{ (x,y) : y \geq 0 \, x \geq -y \}$ as follows:
$$f_H (x,y) := \frac{2}{\pi t^3 (1-t)^3} (x + 2y) e^{-(x+2y)^2 / (2t(1-t))} 
   \, .$$
Then $f_H$ is a conditional density for $(\omega (t) , -M(\omega )(t))$ 
under $\Q_1$ conditioned on $\tau > t$.
\end{lem} 

\noindent{\sc Proof:} Begin with the computation of a density 
for $(\omega , M)$ under the standard Brownian measure $\P$.  
By the reflection principle, using $\P^a$ to denote standard 
Brownian motion starting at $a$, one has 
\begin{eqnarray*}
\P^0 (\omega(t) \in [x , x+dx], M(t) \leq -y) & = & \P^{-2y}
   (\omega (t) \in [x , x+dx]) \\
& = & \sqrt{\frac{1}{2 \pi t}} e^{-(x+2y)^2/(2t)} \, dx \, .
\end{eqnarray*}
Differentiating with respect to $y$ yields
\begin{equation} \label{eq:dens1}
\P^0 (\omega (t) \in [x , x+dx] , M(t) \in [-y , -y + dy]) =
   \frac{4(x+2y)}{2t} \sqrt{\frac{1}{2 \pi t}} e^{-(x+2y)^2 / (2t)}
   \, dx \, dy
\end{equation}
On $R$.

Next, compute the joint density
$$\P (\omega (t) \in [x , x + dx] , M(t) \in [y , y+dy] , \omega
   (1) \in [0 , dz] , \tau > t) \, .$$
To do this, according to the Markov property, one must
multiply~(\ref{eq:dens1}) by
$$\P^x (\omega (1-t) \in [0,dz] , \min_{0 \leq s \leq 1-t} \omega
  (s) \leq -y) \, .$$
By the reflection principle, this last factor is equal to
$$\P^{-2y-x} (\omega (1-t) \in [0,dz]) = \sqrt{\frac{1}{2 \pi t}}
   e^{-(x+2y)^2 / (2(1-t))} \, dz \, ,$$
and multiplying and simplifying $1/t + 1/(1-t)$ to $1/(t(1-t))$ 
in the exponent gives a joint density of
\begin{equation} \label{eq:dens2}
\frac{x + 2y}{\pi \sqrt{t^3(1-t)}} e^{-(x+2y)^2 / (2t(1-t))} \,
   dz \, dx \, dy \, .
\end{equation}

A change of variables simplifies the computation a little.  Let
$u = x+2y$ and $v = (2x - y) / 5$, so that $du \, dv = dx \, dy$
and $x = (u + 10 v)/5, y = (2u - 5v) / 5$.  The region $R$ is
transformed into the region $R' := \{ u \geq 0, (-3/5) u \leq v
\leq (2/5) u$.  Now rewrite the density~(\ref{eq:dens2}) as
\begin{equation} \label{eq:dens22}
\frac{u}{\pi \sqrt{t^3(1-t)}} e^{-u^2 / (2t(1-t))}
   \, dz \, du \, dv \, .
\end{equation}

The conditional density of $(u,v)$ given $\omega (1) \in [0,dz]$
and $\tau < t$, is given by normalizing this.  One must divide by
the integral of~(\ref{eq:dens22}) over $R'$, computed by a simple
linear change of variables $u = r \sqrt{t(1-t)}$ in the third line:
\begin{eqnarray*}
&& \int_0^\infty \int_{-(3/5) u}^{(2/5) u} \frac{u}{\pi \sqrt{t^3(1-t)}}
   e^{-u^2 / (2t(1-t))}  \, dv \, du \, dz \\
& = & \int_0^\infty \frac{u^2}{\pi \sqrt{t^3(1-t)}}
   e^{-u^2 / (2t(1-t))} \, du \, dz \\
& = & \int_0^\infty t(1-t) \frac{r^2}{\pi \sqrt{t^3(1-t)}}
   e^{-r^2/2} \sqrt{t(1-t)} \, dr \, dz \\
& = & \frac{1-t}{\sqrt{2 \pi}} \, dz 
\end{eqnarray*}
using the fact that $\int_0^\infty r^2 e^{-r^2/2} \, dr 
= \sqrt{\pi / 2}$.  Dividing, 
$$ f_H(x(u,v),y(u,v)) = \sqrt{\frac{2}{\pi t^3(1-t)^3}} u
   e^{-u^2 / (2t(1-t))}$$
and plugging in $u = x + 2y$ proves the lemma.   $\Cox$

\noindent{\sc Proof of Theorem~\protect{\ref{th:cond}} continued:}
Let $G$ denote the CDF for the time at which a Brownian bridge on
$[0,1]$ reaches its minimum.  Set the amplitude $2V = 1$ for 
convenience, and note that $\omega (t) - M(t) = x - (-y) = (3u+5v)/5$.
Then by Lemma~\ref{lem:density}, we have
\begin{eqnarray}
\int_\Omega D_1 \, d\Q_1
   & = & \int_\Omega \int_0^1 \left ( \omega (t) - M(t) \right ) 
   \one_{\tau > t} \, dt \, d\Q_1 (\omega ) \nonumber \\
& = & \int_0^1 dt \, (1 - G(t)) \, \int_\Omega \left ( \omega (t) - 
   M (\omega) (t) \right ) \, d(\Q_1 \| \tau > t) (\omega) \nonumber \\
& = & \int_0^1 dt \, (1 - G(t)) \int_{R'} \frac{3u+5v}{5} 
   \, f_H (x(u,v),y(u,v)) \, du \, dv \, . \label{eq:rep1}
\end{eqnarray}

The integral over $R'$ may be computed by substituting $u = r
\sqrt{t(1-t)}$ as before to get
\begin{eqnarray*}
\int_{R'} \frac{3u + 5v}{5} f_H (x,y) dv \, du & = &
   \sqrt{\frac{2}{\pi t^3 (1-t)^3}} \int_0^\infty \, du \, u e^{-u^2 /
   (2t(1-t))} \int_{(-3/5)u}^{(2/5)u} \frac{3u+5v}{5} \, dv \\
& = & \sqrt{\frac{2}{\pi t^3 (1-t)^3}} \int_0^\infty u e^{-u^2 /
   (2t(1-t))} (\frac{3}{5} u^2 - \frac{1}{10} u^2) \, du \\
& = & \sqrt{\frac{1}{2 \pi t^3 (1-t)^3}} \int_0^\infty u^3 e^{-u^2 /
   (2t(1-t))} \, du \\
& = & \sqrt{\frac{1}{2 \pi}} \sqrt{t(1-t)} \int_0^\infty
   e^{-r^2/2} \, dr \\
& = & \frac{1}{2} \sqrt{t(1-t)} \, .
\end{eqnarray*}

Finally, one may evaluate the integral in~(\ref{eq:rep1}) without
knowing the exact distribution function $G$.  By symmetry, we know
$G(t) + G(1-t) = 1$.  The integral over $R'$ is symmetric under
$t \mapsto 1-t$.  Thus
\begin{eqnarray*}
\int_\Omega D_1 \, d\Q_1 & = & \int_0^1 (1 - G(t)) \frac{1}{2}
\sqrt{t(1-t)} \, dt \\
& = & \int_0^1 \frac{1}{4} \sqrt{t(1-t)} \, dt \\
& = & \frac{\pi}{32} \, .
\end{eqnarray*}
Now doubling to get $D$ and multiplying by the amplitude of
$\sqrt{2V}$ leads to
$$\int_\Omega D \, d\Q_{2V} = \frac{\pi \sqrt{2V}}{16} \, ,$$
which establishes the value of $K_2$ and finishes the proof of
Theorem~\ref{th:cond}.

\end{document}